\documentclass{compositio}
\usepackage{mathrsfs}
\usepackage{amsmath}
\usepackage{amssymb}
\usepackage{verbatim}
\usepackage{xcolor}
\usepackage{graphicx}
\usepackage[all]{xy}
\numberwithin{equation}{section}

\theoremstyle{plain}
\newtheorem{theorem}{Theorem}[section]
\newtheorem{proposition}{Proposition}[section]
\newtheorem{corollary}{Corollary}[section]
\newtheorem{lemma}{Lemma}[section]

\theoremstyle{definition}
\newtheorem{definition}{Definition}[section]

\theoremstyle{remark}

\begin{document}

\title{On the relative Nadel-type vanishing theorem}

\author{Jingcao Wu}
\email{wujincao@shufe.edu.cn}
\address{School of Mathematics, Shanghai University of Finance and Economics, Shanghai 200433, People's Republic of China}

\thanks{This research was supported by by NFSC, grant 12271275.}

\begin{abstract}
Let $f:X\rightarrow Y$ be a K\"{a}hler fibration from a complex manifold $X$ to an analytic space $Y$. We show several relative Nadel-type vanishing theorems. 
\end{abstract}

\maketitle

\section{Introduction}
\label{sec:introduction}

The Nadel vanishing theorem \cite{Nad90,Dem93} is a powerful tool in complex geometry. In this paper, we are interested in its relative variant. The start point is the following result. 
\begin{theorem}[(c.f. \cite{Laz04}, Generalizations 9.1.22 and 11.2.15)]\label{t11}
Let $f:X\rightarrow Y$ be a surjective projective morphism of quasi-projective varieties, with $X$ non-singular. Denote by $l$ the dimension of a general fiber of $f$. Let $L$ be a holomorphic line bundle on $X$ with $\kappa(L,f)=l$. Then
\begin{equation*}
R^{q}f_{\ast}(K_{X}\otimes L\otimes\mathscr{I}(f,\|L\|))=0
\end{equation*}
for $q>0$.
\end{theorem}

Here $\mathscr{I}(f,\|L\|)$ is the relative asymptotic multiplier ideal sheaf and the relative Iitaka dimension $\kappa(L,f)$ of $L$ is defined as the maximal integer $m$ such that
\[
\textrm{rank}\,f_{\ast}(L^{k})\geqslant Ck^{m}\quad\textrm{for a universal constant }C\textrm{ and all }k\gg1.
\] 
Note that when $\kappa(L,f)=l$, we will say that $L$ is $f$-big. It is then asked to generalize Theorem \ref{t11} for $L$ possessing weaker positivity. There are several results in this aspect, such as \cite{Mat18}, Theorem 1.2 and \cite{Mat22}, Theorem 1.7. In this paper, we will also provide a relative version of the Nadel-type vanishing theorem. Before stating the main result, we should fix some notations and conventions first. 

Throughout this paper, a fibration $f:X\rightarrow Y$ between analytic spaces refers to a proper, surjective morphism with connected fibers. Moreover, $f$ is said to be K\"{a}hler if $f^{-1}(U)$ is a K\"{a}hler space for any relatively compact, open subset $U$ of $Y$. In the end, we denote by $l$ the dimension of a general fiber of $f$. Then the main theorem is as follows.

\begin{theorem}\label{t12}
Let $f:X\rightarrow Y$ be a K\"{a}hler fibration from a complex manifold $X$ to an analytic space $Y$. Let $(E,h)$ be a smooth Hermitian vector bundle that is Nakano semi-positive, and let $L$ be a holomorphic line bundle with $\kappa(L,f)\geqslant0$ on $X$. Then
\begin{equation*}
R^{q}f_{\ast}(K_{X}\otimes L\otimes E\otimes\mathscr{I}(f,\|L\|))=0
\end{equation*}
for $q>l-\kappa(L,f)$. 
\end{theorem}
Apparently Theorem \ref{t12} extends Theorem \ref{t11}. Next we investigate the higher rank vector bundles equipped with singular metrics, and obtain the following results. 

\begin{theorem}\label{t13}
Let $f:X\rightarrow Y$ be a fibration from a compact K\"{a}hler manifold $X$ to a projective manifold $Y$. Let $(E,h)$ be a singular Hermitian vector bundle that is strong Nakano semi-positive, and let $L$ be an analytic almost base point free line bundle on $X$. Then
\begin{equation*}
R^{q}f_{\ast}(K_{X}\otimes L\otimes\mathcal{E}(h))=0
\end{equation*}
for $q>l-\kappa(L,f)$.   
\end{theorem}

\begin{theorem}\label{t14}
Let $f:X\rightarrow Y$ be a fibration from a projective manifold $X$ to an analytic space $Y$. Let $(E,h)$ be a singular Hermitian vector bundle that is strong Nakano semi-positive, and let $L$ be an $f$-nef line bundle on $X$. Then
\begin{equation*}
R^{q}f_{\ast}(K_{X}\otimes L\otimes\mathcal{E}(h))=0
\end{equation*}
for $q>l-\kappa(L,f)$. 
\end{theorem}

Here $\mathcal{E}(h)$ refers to the multiplier submodule sheaf of $\mathcal{O}(E)$, whose germs are locally $L^{2}$-integrable against $h$. Recall that for a singular Hermitian vector bundle, the description of Nakano semi-positivity \cite{LYZ22} involves the strictly plurisubharmonic functions. It brings geometric obstacle if one attempts to realize them to be the metrics of some line bundles. The strong Nakano semi-positivity is introduced in Sect.\ref{sec:preliminary} in order to get rid of this restriction. These two notions are equivalent at least for the metrics of $C^{2}$. (See Theorem \ref{t21}.) We hope the same thing happens for more singular metrics. $L$ is said to be analytic almost base point free if for any $\varepsilon>0$ and $x\in X$ there exists a singular metric $\varphi$ on $L$ with $i\Theta_{L,\varphi}\geqslant0$ and the Lelong number $\nu(\varphi,x)<\varepsilon$; $L$ is said to be $f$-nef if $L|_{f^{-1}(U)}$ is nef for any Stein open subset $U$ of $Y$. 

In order to prove these theorems, we first apply the Monge--Amp\`{e}re technique developed in \cite{DeP03} and the strong openness property in \cite{LYZ22}, to show the following absolute version.

\begin{theorem}\label{t15}
Let $X$ be a compact K\"{a}hler manifold of dimension $n$.
\begin{enumerate}
    \item[(1)] Assume that $(E,h)$ is a smooth Hermitian vector bundle that is Nakano semi-positive, and $L$ is a holomorphic line bundle with $\kappa(L)\geqslant0$. Then
    \[
    H^{q}(X,K_{X}\otimes L\otimes E\otimes\mathscr{I}(\|L\|))=0
    \]
    for $q>n-\kappa(L)$.
    \item[(2)] Assume that $(E,h)$ is a singular Hermitian vector bundle that is strong Nakano semi-positive, and $L$ is an analytic almost base point free line bundle. Then
    \[
    H^{q}(X,K_{X}\otimes L\otimes\mathcal{E}(h))=0
    \]
    for $q>n-\textrm{nd}(L)$.
\end{enumerate} 
\end{theorem}

Here $\mathscr{I}(\|L\|)$ is the (absolute) asymptotic multiplier ideal sheaf, and $\textrm{nd}(L)$ refers to the numerical dimension, which is uniquely defined when $L$ is nef. So it arises no confusion in our setting. (See Lemma \ref{l21}.) 

Based on Theorem \ref{t15}, we can directly conclude the vanishing result of the higher direct images along the general fiber. Then the torsion-freeness of $R^{q}f_{\ast}(K_{X}\otimes L\otimes E\otimes\mathscr{I}(f,\|L\|))$ implies Theorem \ref{t12}. However, it is not clear whether $R^{q}f_{\ast}(K_{X}\otimes L\otimes\mathcal{E}(h))$ is torsion-free or not. Alternatively we involve some algebraic methods to overcome this problem, which leads to Theorems \ref{t13} and \ref{t14}.

This paper is organized as follows. We first recall some background materials, including the asymptotic multiplier ideal sheaf, the theory of singular Hermitian metrics and so on. In Sect.\ref{sec:vanishing} we prove Theorem \ref{t15}. Then we proceed to show the remaining results in Sect.\ref{sec:main}.

\begin{acknowledgements}
The author would like to thank Prof. Jixiang Fu for the discussion and continuing encouragement.
\end{acknowledgements}

\section{Preliminary}
\label{sec:preliminary}
In this section, unless otherwise stated, we always assume that $X$ is a complex manifold.

\subsection{Asymptotic multiplier ideal sheaf}
This part is mostly collected from \cite{Laz04}. Let $L$ be a holomorphic line bundle on $X$.

First recall the definition of the multiplier ideal sheaf associated to an ideal sheaf $\mathfrak{a}\subset\mathcal{O}_{X}$ and a positive real number $c$. Let $\mu:\tilde{X}\rightarrow X$ be a smooth modification such that $\mu^{\ast}\mathfrak{a}=\mathcal{O}_{\tilde{X}}(-E)$, where $E$ has the simple normal crossing support. Then the multiplier ideal sheaf is defined as
\begin{equation*}
\mathscr{I}(c\cdot\mathfrak{a}):=\mu_{\ast}\mathcal{O}_{\tilde{X}}(K_{\tilde{X}/X}-\lfloor cE\rfloor).
\end{equation*}
Here $K_{\tilde{X}/X}$ is the relative canonical bundle and $\lfloor E\rfloor$ means the round-down.

Now suppose that $\kappa(L)\geqslant0$, and let $\mathfrak{a}_{k}$ be the base-ideal of $L^{k}$. It is easy to verify that for every integer $p\geqslant1$ one has the inclusion
\begin{equation*}
\mathscr{I}(\frac{c}{k}\cdot\mathfrak{a}_{k})\subseteq\mathscr{I}(\frac{c}{pk}\cdot\mathfrak{a}_{pk}).
\end{equation*}
Therefore the family of ideals
\begin{equation*}
\{\mathscr{I}(\frac{c}{k}\cdot\mathfrak{a}_{k})\}_{(k\geqslant0)}
\end{equation*}
has a unique maximal element from the ascending chain condition on ideals.

\begin{definition}
The asymptotic multiplier ideal sheaf associated to $c$ and $L$,
\begin{equation*}
\mathscr{I}(c\|L\|)
\end{equation*}
is defined to be the unique maximal member among the family of ideals $\{\mathscr{I}(\frac{c}{k}\cdot\mathfrak{a}_{k})\}$.
\end{definition}
By definition, there exists some $k$ computing $\mathscr{I}(c\|L\|)$, i.e.
\[
\mathscr{I}(\frac{c}{k}\cdot\mathfrak{a}_{k})=\mathscr{I}(c\|L\|).
\]
Then we can pick $\{u_{1},...,u_{m}\}\subset\Gamma(X,L^{k})$ that generates $\mathfrak{a}_{k}$, and let $\varphi=\frac{1}{k}\log(|u_{1}|^{2}+\cdots+|u_{m}|^{2})$. As a result, 
\[
\mathscr{I}(c\varphi)=\mathscr{I}(c\|L\|).
\]
At this time we call $\varphi$ is associated with $\mathscr{I}(c\|L\|)$. 

The following result will be used later.
\begin{lemma}\label{l21}
We have the following implications:
\[
\mathscr{I}(\|L^{k}\|)=\mathcal{O}_{X}\textrm{ for all }k\Rightarrow L\textrm{ is analytic almost base point free}\Rightarrow L\textrm{ is nef}.
\]
\begin{proof}
The first one is by definition, and the second one is due to \cite{Dem98}.
\end{proof}
\end{lemma}

Next, let $f:X\rightarrow Y$ be a fibration to an analytic space $Y$. Suppose that $\kappa(L,f)\geqslant0$. For a positive integer $k$, there is a naturally defined homomorphism
\begin{equation*}
\rho_{k}:f^{\ast}f_{\ast}(L^{k})\rightarrow L^{k}.
\end{equation*}
The relative base-ideal $\mathfrak{a}_{k,f}$ is then defined as the image of the induced homomorphism
\begin{equation*}
f^{\ast}f_{\ast}L^{k}\otimes L^{-k}\rightarrow\mathcal{O}_{X}.
\end{equation*}
Similarly, we have

\begin{definition}
The relative asymptotic multiplier ideal sheaf associated to $f$, $c$ and $L$,
\begin{equation*}
\mathscr{I}(f,c\|L\|)
\end{equation*}
is defined to be the unique maximal member among the family of ideals $\{\mathscr{I}(\frac{c}{k}\cdot\mathfrak{a}_{k,f})\}$.
\end{definition}

\subsection{Singular Hermitian metrics}
Let $p:E\rightarrow X$ be a holomorphic vector bundle of rank $r$ over $X$.

\begin{definition}(\cite{BP08,PaT18,Rau15})
A singular Hermitian metric on $E$ is a function $h$ that associates to every point $x\in X$ a singular Hermitian inner product $|\cdot|_{h,x}:E_{x}\rightarrow[0,+\infty]$ on the complex vector space $E_{x}$, subject to the following two conditions:
\begin{enumerate}
\item[1.] $h$ is finite and positive definite almost everywhere, meaning that for all $x$ outside a set of Lebesgue measure zero, $|\cdot|_{h,x}$ is a Hermitian inner product on $E_{x}$;
\item[2.] $h$ is measurable, meaning that the function
\[
\begin{split}
|F|_{h}:U&\rightarrow[0,+\infty] \\
x&\mapsto|F(x)|_{h,x}
\end{split}
\]
is measurable for any open $U\subseteq X$ and $F\in\Gamma(U,E)$.
\end{enumerate}
\end{definition}

\begin{definition}(\cite{Rau15})
Let $h$ be a singular Hermitian metric on $E$, which canonically induces a singular metric $h^{\ast}$ on the dual bundle $E^{\ast}$. 
\begin{enumerate}
\item[(1)] $(E,h)$ is called Griffiths semi-negative (or negatively curved) if for any (local) section $F$ of $E$, the function $\log|F|^{2}_{h}$ is plurisubharmonic.
\item[(2)] $(E,h)$ is called Griffiths semi-positive (or positively curved) if $(E^{\ast},h^{\ast})$ is Griffiths semi-negative. 
\end{enumerate}
\end{definition}

The definition for the Nakano-type positivity would be subtle. Firstly we introduce the following conventions. Fix a Hermitian metric $\omega$ and a Hermitian line bundle $(L,\phi)$ on $X$. Denote by $\langle\cdot,\cdot\rangle_{\phi,h,\omega}$ (resp. $|\cdot|_{\phi,h,\omega}$) the pointwise inner product (resp. norm) on $(L\otimes E)$-valued $(p,q)$-forms induced by $\phi$, $h$ and $\omega$, and by $(\cdot,\cdot)_{\phi,h,\omega}:=\int_{X}\langle\cdot,\cdot\rangle_{\phi,h,\omega}dV_{\omega}$ (resp. $\|\cdot\|_{\phi,h,\omega}:=(\int_{X}|\cdot|^{2}_{\phi,h,\omega}dV_{\omega})^{1/2}$) the $L^{2}$-inner product (resp. $L^{2}$-norm). Let $L^{p,q}_{(2)}(X,L\otimes E)_{\phi,h,\omega}$ be the collection of $(L\otimes E)$-valued $(p,q)$-forms with measurable coefficients and bounded $L^{2}$-norm. 

\begin{definition}(\cite{DNWZ23})
A singular Hermitian vector bundle $(E,h)$ is said to satisfy the optimal $L^{2}$-estimate condition, if for any positive line bundle $(A,\phi)$, and any compactly supported, $\bar{\partial}$-closed $\beta\in L^{n,1}_{(2)}(X,A\otimes E)_{\phi,h,\omega}$ with
\[
(B^{-1}_{\phi,\omega}\beta,\beta)_{\phi,h,\omega}<\infty,
\]
there exists $u\in L^{n,0}_{(2)}(X,A\otimes E)_{\phi,h,\omega}$ satisfying $\beta=\bar{\partial}u$ and
\[
\|u\|^{2}_{\phi,h,\omega}\leqslant(B^{-1}_{\phi,\omega}\beta,\beta)_{\phi,h,\omega}.
\]
Here $B_{\phi,\omega}=[i\partial\bar{\partial}\phi\otimes\mathrm{Id}_{E},\Lambda_{\omega}]$.
\end{definition}

\begin{definition}(\cite{LYZ22})\label{d26}
$(E,h)$ is called Nakano semi-positive if
\begin{enumerate}
    \item[1.] $(E,h)$ is Griffiths semi-positive;
    \item[2.] $(E,h)$ satisfies the optimal $L^{2}$-estimate condition.
\end{enumerate}
\end{definition}

In the same paper, they showed that when $h$ is at least of $C^{2}$, it coincides with the classic Nakano semi-positivity described via curvature. From this point of view, Definition \ref{d26} is a good attempt to generalize the Nakano-type positivity into the singular setting. However, it would be a little bit luxurious in practice to involve a positive line bundle. So we make the following adjustment.
 
\begin{definition}
A singular Hermitian vector bundle $(E,h)$ is said to satisfy the strong $L^{2}$-estimate condition, if for any collection of the following data: 
\begin{enumerate}
    \item[1.] A set of positive real numbers $\{\varepsilon\}$ that decreases to zero;
    \item[2.] A holomorphic line bundle $L$ equipped with a family of (singular) metrics $\{\phi_{\varepsilon}\}$, whose curvature satisfying $i\Theta_{L,\phi_{\varepsilon}}\geqslant-\varepsilon\omega$;
    \item[3.] $\lambda\geqslant0$ that makes $\Theta:=i\Theta_{L,\phi_{\varepsilon}}+\lambda\omega$ a strictly positive $(1,1)$-current;
    \item[4.] $\bar{\partial}$-closed $\beta\in L^{n,q}_{(2)}(X,L\otimes E)_{\phi_{\varepsilon},h,\omega}$ with
\[
(B^{-1}_{\Theta,\omega}\beta,\beta)_{\phi_{\varepsilon},h,\omega}<\infty,
\] 
\end{enumerate}
there exist elements $u_{\varepsilon},v_{\varepsilon}$ and positive constant $C$ (independent of $\varepsilon$) satisfying $\beta=\bar{\partial}u_{\varepsilon}+v_{\varepsilon}$ and
\[
\|u_{\varepsilon}\|^{2}_{\phi_{\varepsilon},h,\omega}+\frac{C}{\varepsilon}\|v_{\varepsilon}\|^{2}_{\phi_{\varepsilon},h,\omega}\leqslant(B^{-1}_{\Theta,\omega}\beta,\beta)_{\phi_{\varepsilon},h,\omega}.
\]
Here $B_{\Theta,\omega}=[\Theta\otimes\mathrm{Id}_{E},\Lambda_{\omega}]$.
\end{definition}

Note that for a positive line bundle $(A,\phi)$, if exists, we can take $\Theta=i\partial\bar{\partial}\phi>\delta\omega$ for a positive $\delta$, then for any $(n,1)$-form $\beta$ we have
\[
\langle B^{-1}_{\phi,\omega}\beta,\beta\rangle_{\phi,h,\omega}\leqslant\frac{1}{a_{1}}|\beta|^{2}_{\phi,h,\omega}<\frac{1}{\delta}|\beta|^{2}_{\phi,h,\omega}.
\]
Here $\delta<a_{1}\leqslant\cdots\leqslant a_{n}$ are the eigenvalues of $\Theta$ at every point $x\in X$ with respect to $\omega(x)$. Hence
\[
\|v_{\varepsilon}\|^{2}_{\phi,h,\omega}\leqslant\frac{\varepsilon}{C\delta}\|\beta\|^{2}_{\phi,h,\omega},
\]
which converges to zero as $\varepsilon$ tends to zero. Therefore the strong $L^{2}$-estimate condition implies the optimal $L^{2}$-estimate condition. Moreover, we can show that 

\begin{theorem}\label{t21}
Let $(X,\omega)$ be a complete K\"{a}hler manifold. If $h$ is of $C^{2}$, then the following statements are equivalent:
\begin{enumerate}
    \item[(1)] $(E,h)$ satisfies the strong $L^{2}$-estimate condition.
    \item[(2)] $(E,h)$ satisfies the optimal $L^{2}$-estimate condition.
    \item[(3)] $(E,h)$ is Nakano semi-positive (in the classic sense). 
\end{enumerate}
\begin{proof}
The equivalence between (2) and (3) is established in \cite{DNWZ23}. If $(E,h)$ is Nakano semi-positive, the standard $L^{2}$-estimate implies that $(E,h)$ satisfies the strong $L^{2}$-estimate condition. One could also consult our proof of Theorem \ref{t15} for a suggested argument. It remains to prove (1) $\Rightarrow$ (3).

Indeed, fix a set of positive real numbers $\{\varepsilon\}$ that decreases to zero. Let $L$ be the trivial line bundle, and let $\phi_{\varepsilon}$ be the natural flat metric on $L$ for every $\varepsilon$. Let $\Theta=\omega$. By the strong $L^{2}$-estimate condition for an arbitrary compactly supported, $\bar{\partial}$-closed $\beta\in L^{n,q}_{(2)}(X,E)_{h,\omega}$ with
\[
(B^{-1}_{\Theta,\omega}\beta,\beta)_{h,\omega}<\infty,
\]
there exist elements $u_{\varepsilon},v_{\varepsilon}$ and positive constant $C$ satisfying $\beta=\bar{\partial}u_{\varepsilon}+v_{\varepsilon}$ and
\[
\|u_{\varepsilon}\|^{2}_{h,\omega}+\frac{C}{\varepsilon}\|v_{\varepsilon}\|^{2}_{h,\omega}\leqslant(B^{-1}_{\Theta,\omega}\beta,\beta)_{h,\omega}.
\] 
Now for any test $E$-valued $(n,q)$-form $\alpha$, we have
\[
\begin{split}
|(\alpha,\beta)_{h,\omega}|&=|(\alpha,\bar{\partial}u_{\varepsilon}+v_{\varepsilon})_{h,\omega}| \\
&\leqslant|(\bar{\partial}^{\ast}\alpha,u_{\varepsilon})_{h,\omega}|+|(\alpha,v_{\varepsilon})_{h,\omega}| \\
&\leqslant\|\bar{\partial}^{\ast}\alpha\|_{h,\omega}\|u_{\varepsilon}\|_{h,\omega}+|(\alpha,v_{\varepsilon})_{h,\omega}|.
\end{split}
\]
Whereas 
\begin{equation}\label{e21}
\begin{split}
\|\bar{\partial}^{\ast}\alpha\|^{2}_{h,\omega}\|u_{\varepsilon}\|^{2}_{h,\omega}\leqslant&(\|D^{\ast}_{1}\alpha\|^{2}_{h,\omega}+([i\Theta_{E,h},\Lambda_{\omega}]\alpha,\alpha)_{h,\omega}-\|\bar{\partial}\alpha\|^{2}_{h,\omega})\cdot(B^{-1}_{\Theta,\omega}\beta,\beta)_{h,\omega},
\end{split}
\end{equation}
and
\begin{equation}\label{e22}
\begin{split}
\lim_{\varepsilon\rightarrow0}|(\alpha,v_{\varepsilon})_{h,\omega}|^{2}&\leqslant\lim_{\varepsilon\rightarrow0}\|\alpha\|^{2}_{h,\omega}\|v_{\varepsilon}\|^{2}_{h,\omega} \\
&\leqslant\lim_{\varepsilon\rightarrow0}\frac{\varepsilon}{C}\|\alpha\|^{2}_{h,\omega}\|\beta\|^{2}_{h,\omega}=0.
\end{split}
\end{equation}
Here $D^{\ast}_{1}$ refers to the adjoint of the $(1,0)$-part of the Chern connection. Hence (\ref{e21}) and (\ref{e22}) together imply that 
\begin{equation}\label{e23}
\begin{split}
|(\alpha,\beta)_{h,\omega}|^{2}\leqslant&(\|D^{\ast}_{1}\alpha\|^{2}_{h,\omega}+([i\Theta_{E,h},\Lambda_{\omega}]\alpha,\alpha)_{h,\omega}-\|\bar{\partial}\alpha\|^{2}_{h,\omega})\cdot( B^{-1}_{\Theta,\omega}\beta,\beta)_{h,\omega}.
\end{split}
\end{equation}
Then let $\alpha=B^{-1}_{\Theta,\omega}\beta$, (\ref{e23}) becomes
\begin{equation}\label{e24}
\begin{split}
0&\leqslant( B_{\Theta,\omega}\alpha,\alpha)_{h,\omega}\leqslant(\|D^{\ast}_{1}\alpha\|^{2}_{h,\omega}+([i\Theta_{E,h},\Lambda_{\omega}]\alpha,\alpha)_{h,\omega}-\|\bar{\partial}\alpha\|^{2}_{h,\omega}).
\end{split}
\end{equation}
Suppose that $(E,h)$ is not Nakano semi-positive. The same argument as \cite{DNWZ23}, Theorem 1.1 then leads to a contradiction to (\ref{e24}). The proof is complete.
\end{proof}
\end{theorem}
It would be valuable to ask that for a general singular metric, whether the optimal $L^{2}$-estimate condition implies the strong $L^{2}$-estimate or not, which will make the following definition equivalent to Definition \ref{d26}.

\begin{definition}
$(E,h)$ is called strong Nakano semi-positive if
\begin{enumerate}
    \item[1.] $(E,h)$ is Griffiths semi-positive;
    \item[2.] $(E,h)$ satisfies the strong $L^{2}$-estimate condition.
\end{enumerate}
\end{definition}

Then we collect several properties concerning the singular Hermitian metrics, which will be used later. 

\begin{proposition}[(\cite{Ina22})]
If $(E,h)$ is Nakano semi-positive, then $\mathcal{E}(h)$ is coherent.
\end{proposition}

\begin{proposition}[(\cite{Rau15}, Proposition 1.3)]\label{p22}
If $(E,h)$ is Griffiths semi-negative, then $\log\det h$ is a plurisubharmonic function.
\end{proposition}

\begin{proposition}[(\cite{LYZ22}, Lemma 4.1)]\label{p23}
If $(E,h)$ is Nakano semi-positive, then for a given section $F$ of $E$, the following are equivalent:
\begin{enumerate}
    \item[1.] $|F|^{2}_{h(\det h)^{s}}$ is locally integrable for some $s>0$;
    \item[2.] $|F|^{2t}_{h}$ is locally integrable for some $t>1$.
\end{enumerate}
\end{proposition}

\begin{proposition}[(\cite{LYZ22}, Corollary 4.4)]\label{p24}
Let $\psi$ be a plurisubharmonic function. If $(E,h)$ is Nakano semi-positive, then
\[
\cup_{s>0}\mathcal{E}(he^{-s\psi})=\mathcal{E}(h).
\]
\end{proposition}

In particular, combining with Propositions \ref{p22}, \ref{p23} and \ref{p24}, we conclude
\begin{corollary}\label{c21}
If $(E,h)$ is Nakano semi-positive, and $F$ is a section of $E$ such that $|F|^{2}_{h}$ is locally integrable, then there exists some $t>1$ such that $|F|^{2t}_{h}$ is locally integrable.
\end{corollary}

\section{A Nadel-type vanishing theorem}
\label{sec:vanishing}

The following is a variant of \cite{Wu22}, Theorem 1.3.

\begin{theorem}[(=Theorem \ref{t15})]\label{t31}
Let $X$ be a compact K\"{a}hler manifold $X$ of dimension $n$.
\begin{enumerate}
    \item[(1)] Assume that $(E,h)$ is a smooth Hermitian vector bundle that is Nakano semi-positive, and $L$ is a holomorphic line bundle with $\kappa(L)\geqslant0$. Then
    \[
    H^{q}(X,K_{X}\otimes L\otimes E\otimes\mathscr{I}(\|L\|))=0
    \]
    for $q>n-\kappa(L)$.
    \item[(2)] Assume that $(E,h)$ is a singular Hermitian vector bundle that is strong Nakano semi-positive, and $L$ is analytic almost base point free. Then
    \[
    H^{q}(X,K_{X}\otimes L\otimes\mathcal{E}(h))=0
    \]
    for $q>n-\textrm{nd}(L)$.
\end{enumerate} 
\begin{proof}
(1) We firstly do some reductions. Let $\varphi$ be a metric on $L$ that is associated with $\mathscr{I}(\|L\|)$. Then $\varphi$ is of the form that
\[
\varphi=\frac{1}{p}\log(|u_{1}|^{2}+\cdots+|u_{m}|^{2}),
\]
where $\{u_{1},...,u_{m}\}$ is a basis of $H^{0}(X,L^{p})$ for some $p\gg0$ and divisible enough. Let $\mathfrak{a}_{p}$ be the base-ideal of the linear series $|L^{p}|$, and let $\mu:\tilde{X}\rightarrow X$ be the log-resolution of $|L^{p}|$ such that $\mu^{\ast}|L^{p}|=|V_{p}|+F_{p}$. Accordingly,
\[
\mu^{\ast}\mathfrak{a}_{p}=\mathcal{O}_{\tilde{X}}(-F_{p}),
\]
where $F_{p}=\sum\lambda_{i}G_{i}$ has simple normal crossing support. Then locally
\[
\varphi\circ\mu=\frac{1}{p}\sum\lambda_{i}\log|g_{i}|^{2}+\tau,
\]
where $g_{i}$ is the local generator of $G_{i}$ and $\tau$ is a (local) smooth quasi-plurisubharmonic function. Hence 
\[
\mathscr{I}(\varphi)=\mu_{\ast}\mathcal{O}_{\tilde{X}}(K_{\tilde{X}/X}-\sum\lfloor\frac{\lambda_{i}}{p}\rfloor G_{i}).
\]
Now, a direct application of the Leray spectral sequence \cite{Har77} implies that
\[
H^{q}(X,K_{X}\otimes L\otimes E\otimes\mathscr{I}(\varphi))=H^{q}(\tilde{X},K_{\tilde{X}}\otimes\mu^{\ast}L\otimes\mu^{\ast}E\otimes\mathcal{O}_{\tilde{X}}(-\sum\lfloor\frac{\lambda_{i}}{p}\rfloor G_{i})).
\] 
Let $\hat{L}:=\mu^{\ast}L\otimes\mathcal{O}_{\tilde{X}}(-\sum\lfloor\frac{\lambda_{i}}{p}\rfloor G_{i})$. Then it is enough to consider
\[
H^{q}(\tilde{X},K_{\tilde{X}}\otimes\hat{L}\otimes\mu^{\ast}E).
\]

Let $\tilde{L}=\mu^{\ast}L\otimes\mathcal{O}_{\tilde{X}}(-\sum\frac{\lambda_{i}}{p} G_{i})$, which is a $\mathbb{Q}$-bundle. Then 
\[
\hat{L}=\tilde{L}\otimes\mathcal{O}_{\tilde{X}}(\sum\{\frac{\lambda_{i}}{p}\}G_{i}),
\]
where $\{\frac{\lambda_{i}}{p}\}$ refers to the fractional part of $\frac{\lambda_{i}}{p}$. Observe that $\tilde{L}^{p}$, which is a $\mathbb{Z}$-bundle, is semi-positive equipped with the smooth metric $\psi=p\varphi\circ\mu-\sum\lambda_{i}\log|g_{i}|^{2}$. 

On the other hand, we have
\[
\mu^{\ast}L^{p}=\tilde{L}^{p}\otimes\mathcal{O}_{\tilde{X}}(F_{p}).
\]
Then the morphism $\phi_{|\tilde{L}^{p}|}:\tilde{X}\rightarrow\mathbb{P}^{N}$ defined by the free linear series $|\tilde{L}^{p}|$ resolves the indeterminacy of $\phi_{|L^{p}|}$, and in particular it is generically finite over its image. Therefore 
\[
\kappa(\tilde{L}^{p})=\kappa(L^{p})=\kappa(L).
\] 
As a result, the numerical dimension $\textrm{nd}(\tilde{L}^{p})$, which is uniquely defined since $\tilde{L}^{p}$ is nef, is not less than $\kappa(L)$. Hence $c_{1}(\tilde{L}^{p})^{\kappa(L)}\neq0$. 

Fix a K\"{a}hler metric $\omega$ on $\tilde{X}$. Let $\varepsilon>0$. Then $i\Theta_{\tilde{L}^{p},\psi}+\varepsilon\omega$ is a K\"{a}hler form, hence by the Calabi--Yau theorem for complex Monge--Amp\`{e}re equations there exists a smooth metric $\varphi_{\varepsilon}$ on $\tilde{L}^{p}$ such that 
\begin{equation}\label{e31}
i\Theta_{\tilde{L}^{p},\varphi_{\varepsilon}}+\varepsilon\omega>0\quad\textrm{and}\quad(i\Theta_{\tilde{L}^{p},\varphi_{\varepsilon}}+\varepsilon\omega)^{n}=C_{\varepsilon}\omega^{n}.
\end{equation}
Here $C_{\varepsilon}>0$ is a constant such that
\[
C_{\varepsilon}=\frac{\int_{\tilde{X}}(c_{1}(\tilde{L}^{p})+\varepsilon\omega)^{n}}{\int_{\tilde{X}}\omega^{n}}\geqslant C\varepsilon^{n-\kappa(L)}
\]
for a universal constant $C$. Now we consider the metric
\[
\phi_{\varepsilon}=\frac{1}{p}(\delta\varphi_{\varepsilon}+(1-\delta)\psi)+\sum\{\frac{\lambda_{i}}{p}\}\log|g_{i}|^{2}
\]
on $\hat{L}$, where $\delta>0$ is a sufficiently small number which will be fixed later. Since $\varphi_{\varepsilon}$ and $\psi$ are both smooth and $\{\frac{\lambda_{i}}{p}\}<1$, the $L^{2}$-norm $\|\alpha\|^{2}_{\phi_{\varepsilon}}$ is always bounded for any $\alpha\in A^{s,t}(\tilde{X},\hat{L})$. Moreover, by construction,
\[
\begin{split}
i\Theta_{\hat{L},\phi_{\varepsilon}}+\frac{2\varepsilon}{p}\omega&=\frac{1}{p}(\delta(i\Theta_{\tilde{L}^{p},\varphi_{\varepsilon}}+\varepsilon\omega)+(1-\delta)(i\Theta_{\tilde{L}^{p},\psi}+\varepsilon\omega))+\sum\{\frac{\lambda_{i}}{p}\}[G_{i}]+\frac{\varepsilon}{p}\omega \\ 
&\geqslant\frac{\delta}{p}(i\Theta_{\tilde{L}^{p},\varphi_{\varepsilon}}+\varepsilon\omega)+\frac{\varepsilon}{p}\omega.
\end{split}
\]
Here $[G_{i}]$ is the current of integration.

Denote by $0<a_{1}\leqslant\cdots\leqslant a_{n}$ and $0<\hat{a}_{1}\leqslant\cdots\leqslant\hat{a}_{n}$, respectively, the eigenvalues of the curvature forms $i\Theta_{\tilde{L}^{p},\varphi_{\varepsilon}}+\varepsilon\omega$ and $i\Theta_{\hat{L},\phi_{\varepsilon}}+\frac{2\varepsilon}{p}\omega$ at every point $x\in\tilde{X}$, with respect to the base K\"{a}hler metric $\omega(x)$. We find $\hat{a}_{j}\geqslant\frac{\delta}{p}a_{j}+\frac{\varepsilon}{p}$. On the other hand the Monge--Amp\`{e}re equation (\ref{e31}) tells us that
\begin{equation}\label{e32}
a_{1}\cdots a_{n}=C_{\varepsilon}\geqslant C\varepsilon^{n-\kappa(L)}
\end{equation}
everywhere on $\tilde{X}$. 

We apply the basic Bochner--Kodaira inequality for every $\hat{L}\otimes\mu^{\ast}E$-valued $(n,q)$-form to obtain
\begin{equation}\label{e33}
\begin{split}
\|\bar{\partial}\alpha\|^{2}_{\phi_{\varepsilon},\mu^{\ast}h,\omega}+\|\bar{\partial}^{\ast}\alpha\|^{2}_{\phi_{\varepsilon},\mu^{\ast}h,\omega}&\geqslant(i\Theta_{\hat{L}\otimes\mu^{\ast}E,\phi_{\varepsilon},\mu^{\ast}h}\Lambda\alpha,\alpha)_{\phi_{\varepsilon},\mu^{\ast}h,\omega} \\
&\geqslant\int_{\tilde{X}}(\hat{a}_{1}+\cdots+\hat{a}_{q}-\frac{2q\varepsilon}{p})|\alpha|^{2}_{\phi_{\varepsilon},\mu^{\ast}h,\omega}dV_{\omega}.
\end{split}
\end{equation}
Note that although $\phi_{\varepsilon}$ is not smooth on the whole $\tilde{X}$, it is explained in \cite{DeP03} that in the limit we can recover the same estimates as if we were in the smooth case.

Now we are ready to prove the desired vanishing result, namely for any $q>n-\kappa(L)$,
\[
H^{q}(X,K_{X}\otimes L\otimes E\otimes\mathscr{I}(\|L\|))=H^{q}(\tilde{X},K_{\tilde{X}}\otimes\hat{L}\otimes\mu^{\ast}E)=0.
\]
Let us take a cohomology class $[\beta]\in H^{q}(\tilde{X},K_{\tilde{X}}\otimes\hat{L}\otimes\mu^{\ast}E)$. By using the de Rham--Weil isomorphism, we obtain a representative $\beta$ which is a smooth $\hat{L}\otimes\mu^{\ast}E$-valued $(n,q)$-form. We want to show that $\beta$ is actually a boundary. Since $H^{q}(\tilde{X},K_{\tilde{X}}\otimes\hat{L}\otimes\mu^{\ast}E)$ is a finite dimensional Hausdorff vector space whose topology is induced by the $L^{2}$-Hilbert space topology on the space of forms, it is enough to show that we can approach $\beta$ by $\bar{\partial}$-exact forms in $\|\cdot\|_{\infty,\mu^{\ast}h,\omega}$. Certainly $\|\cdot\|_{\infty,\mu^{\ast}h,\omega}$ here means the $L^{2}$-norm defined by a reference metric $\phi_{\infty}$ on $\hat{L}$, which implicitly fixed at the beginning.

As in H\"{o}rmander \cite{Hor65}, we write every form $\alpha$ in the domain of the $L^{2}$-extension (with respect to $\|\cdot\|_{\phi_{\varepsilon},\mu^{\ast}h,\omega}$) of $\bar{\partial}^{\ast}$ as $\alpha=\alpha_{1}+\alpha_{2}$ with
\[
\alpha_{1}\in\textrm{Ker}\bar{\partial}\quad\textrm{and}\quad\alpha_{2}\in(\textrm{Ker}\bar{\partial})^{\perp}=\overline{\textrm{Im}\bar{\partial}^{\ast}}\subset\textrm{Ker}\bar{\partial}^{\ast}.
\] 
Therefore, since $\beta\in\textrm{Ker}\bar{\partial}$,
\[
\begin{split}
|(\beta,\alpha)_{\phi_{\varepsilon},\mu^{\ast}h,\omega}|^{2}&=|(\beta,\alpha_{1})_{\phi_{\varepsilon},\mu^{\ast}h,\omega}|^{2} \\
&\leqslant(\int_{\tilde{X}}\frac{1}{\hat{a}_{1}+\cdots+\hat{a}_{q}}|\beta|^{2}_{\phi_{\varepsilon},\mu^{\ast}h,\omega}dV_{\omega})(\int_{\tilde{X}}(\hat{a}_{1}+\cdots+\hat{a}_{q})|\alpha_{1}|^{2}_{\phi_{\varepsilon},\mu^{\ast}h,\omega}dV_{\omega}).
\end{split}
\]
As $\bar{\partial}\alpha_{1}=0$, an application of (\ref{e33}) to $\alpha_{1}$ shows that the second integral in the right-hand side is bounded above by
\[
\|\bar{\partial}^{\ast}\alpha_{1}\|^{2}_{\phi_{\varepsilon},\mu^{\ast}h,\omega}+\frac{2q\varepsilon}{p}\|\alpha_{1}\|^{2}_{\phi_{\varepsilon},\mu^{\ast}h,\omega}\leqslant\|\bar{\partial}^{\ast}\alpha\|^{2}_{\phi_{\varepsilon},\mu^{\ast}h,\omega}+\frac{2q\varepsilon}{p}\|\alpha\|^{2}_{\phi_{\varepsilon},\mu^{\ast}h,\omega},
\]
so we finally get
\[
|(\beta,\alpha)_{\phi_{\varepsilon},\mu^{\ast}h,\omega}|^{2}\leqslant(\int_{\tilde{X}}\frac{1}{\hat{a}_{1}+\cdots+\hat{a}_{q}}|\beta|^{2}_{\phi_{\varepsilon},\mu^{\ast}h,\omega}dV_{\omega})(\|\bar{\partial}^{\ast}\alpha\|^{2}_{\phi_{\varepsilon},\mu^{\ast}h,\omega}+\frac{2q\varepsilon}{p}\|\alpha\|^{2}_{\phi_{\varepsilon},\mu^{\ast}h,\omega}).
\]
By the Hahn--Banach theorem, we can find elements $u_{\varepsilon}$, $v_{\varepsilon}$ such that
\[
(\beta,\alpha)_{\phi_{\varepsilon},\mu^{\ast}h,\omega}=(u_{\varepsilon},\bar{\partial}^{\ast}\alpha)_{\phi_{\varepsilon},\mu^{\ast}h,\omega}+(v_{\varepsilon},\alpha)_{\phi_{\varepsilon},\mu^{\ast}h,\omega}
\]
for all $\alpha$. Namely $\beta=\bar{\partial}u_{\varepsilon}+v_{\varepsilon}$, with
\[
\|u_{\varepsilon}\|^{2}_{\phi_{\varepsilon},\mu^{\ast}h,\omega}+\frac{p}{2q\varepsilon}\|v_{\varepsilon}\|^{2}_{\phi_{\varepsilon},\mu^{\ast}h,\omega}\leqslant\int_{X}\frac{1}{\hat{a}_{1}+\cdots+\hat{a}_{q}}|\beta|^{2}_{\phi_{\varepsilon},\mu^{\ast}h,\omega}dV_{\omega}.
\]
As a consequence, the $L^{2}$-distance of $\beta$ to the space of $\bar{\partial}$-exact forms is bounded by $\|v_{\varepsilon}\|_{\phi_{\varepsilon},\mu^{\ast}h,\omega}$, where
\[
\|v_{\varepsilon}\|^{2}_{\phi_{\varepsilon},\mu^{\ast}h,\omega}\leqslant\frac{2q\varepsilon}{p}\int_{\tilde{X}}\frac{1}{\hat{a}_{1}+\cdots+\hat{a}_{q}}|\beta|^{2}_{\phi_{\varepsilon},\mu^{\ast}h,\omega}dV_{\omega}.
\]

We normalize the choice of the potentials $\phi_{\varepsilon}$ so that $\sup_{X}\phi_{\varepsilon}=0$. From this we infer that
\[
\int_{\tilde{X}}|v_{\varepsilon}|^{2}_{\infty,\mu^{\ast}h,\omega}dV_{\omega}\leqslant\int_{\tilde{X}}|v_{\varepsilon}|^{2}_{\phi_{\varepsilon},\mu^{\ast}h,\omega}dV_{\omega}\leqslant\frac{2q\varepsilon}{p}\int_{\tilde{X}}\frac{1}{\hat{a}_{1}+\cdots+\hat{a}_{q}}|\beta|^{2}_{\phi_{\varepsilon},\mu^{\ast}h,\omega}dV_{\omega}.
\]
It remains to show that the right-hand side converges to zero. By construction $\hat{a}_{j}\geqslant\frac{\delta}{p}a_{j}+\frac{\varepsilon}{p}$ and (\ref{e32}) implies
\[
a^{q}_{q}a_{q+1}\cdots a_{n}\geqslant a_{1}\cdots a_{n}\geqslant C\varepsilon^{n-\kappa(L)},
\] 
hence
\[
\frac{1}{a_{1}+\cdots+a_{q}}\leqslant\frac{1}{a_{q}}\leqslant C^{-1/q}\varepsilon^{-(n-\kappa(L))/q}(a_{q+1}\cdots a_{n})^{1/q}.
\]
We infer
\[
\gamma_{\varepsilon}:=\frac{q\varepsilon}{p(\hat{a}_{1}+\cdots+\hat{a}_{q})}\leqslant\min(1,\frac{q\varepsilon}{\delta a_{q}})\leqslant\min(1,C^{\prime}\delta^{-1}\varepsilon^{1-(n-\kappa(L))/q}(a_{q+1}\cdots a_{n})^{1/q}).
\]
Notice that
\[
\int_{\tilde{X}}a_{q+1}\cdots a_{n}dV_{\omega}\leqslant\int_{\tilde{X}}(i\Theta_{\tilde{L}^{p},\varphi_{\varepsilon}}+\varepsilon\omega)^{n-q}\wedge\omega^{q}=(c_{1}(\tilde{L}^{p})+\varepsilon[\omega])^{n-q}[\omega]^{q}\leqslant C^{\prime\prime},
\]
hence the functions $(a_{q+1}\cdots a_{n})^{1/q}$ are uniformly bounded in $L^{1}$-norm as $\varepsilon$ tens to zero. Since $1-\frac{n-\kappa(L)}{q}>0$ by hypothesis, we conclude that $\gamma_{\varepsilon}$ converges almost everywhere to zero as $\varepsilon$ tends to zero.

Now
\[
\begin{split}
\frac{2q\varepsilon}{p}\int_{\tilde{X}}\frac{1}{\hat{a}_{1}+\cdots+\hat{a}_{q}}|\beta|^{2}_{\phi_{\varepsilon},\mu^{\ast}h,\omega}dV_{\omega}&=2\int_{\tilde{X}}\gamma_{\varepsilon}|\beta|^{2}_{\infty,\mu^{\ast}h,\omega}e^{-\frac{1}{p}(\delta\varphi_{\varepsilon}+(1-\delta)\psi)}\frac{1}{\Pi|g_{i}|^{2\{\frac{\lambda_{i}}{p}\}}}dV_{\omega} \\
&\leqslant2(\int_{\tilde{X}}e^{-\frac{s\delta}{p}\varphi_{\varepsilon}}dV_{\omega})^{1/s}(\int_{\tilde{X}}\frac{\gamma^{t}_{\varepsilon}|\beta|^{2t}_{\infty,\mu^{\ast}h,\omega}e^{-\frac{t(1-\delta)}{p}\psi}}{\Pi|g_{i}|^{2t\{\frac{\lambda_{i}}{p}\}}}dV_{\omega})^{1/t}
\end{split}
\]
for any $s,t>1$ with $\frac{1}{s}+\frac{1}{t}=1$ by H\"{o}lder's inequality. Then we can take $t=1+\xi$ with $\xi$ small enough such that
\[
\int_{\tilde{X}}\frac{1}{\Pi|g_{i}|^{2t\{\frac{\lambda_{i}}{p}\}}}dV_{\omega}<+\infty.
\]
As $\gamma_{\varepsilon}\leqslant1$, the Lebesgue dominated convergence theorem shows that for every $1>\delta>0$,
\[
\int_{\tilde{X}}\frac{\gamma^{t}_{\varepsilon}|\beta|^{2t}_{\infty,\mu^{\ast}h,\omega}e^{-\frac{t(1-\delta)}{p}\psi}}{\Pi|g_{i}|^{2t\{\frac{\lambda_{i}}{p}\}}}dV_{\omega}
\]
converges to zero as $\varepsilon$ tends to zero. On the other hand, since the curvature forms
\[
i\Theta_{\tilde{L}^{p},\varphi_{\varepsilon}}>-\varepsilon\omega
\]
all sit in $c_{1}(\tilde{L}^{p})$, the family of quasi-plurisubharmonic functions $\{\varphi_{\varepsilon}\}$ is a bounded family (after normalization) with respect to the $L^{1}$-norm. By standard results of Skoda \cite{Sko72}, we conclude that there exists a small constant $\eta>0$ such that $\int_{\tilde{X}}e^{-\eta\varphi_{\varepsilon}}$ is uniformly bounded. By choosing $\delta\leqslant\frac{p\eta}{s}$, the integral $\int_{\tilde{X}}e^{-\frac{s\delta}{p}\varphi_{\varepsilon}}$ remains bounded. The proof is then complete.  

(2) Fix a K\"{a}hler metric $\omega$ on $X$. Let $\varepsilon>0$ be an arbitrary real number. Since $L$ is nef by Lemma \ref{l21}, $c_{1}(L)+\varepsilon[\omega]$ is a K\"{a}hler class. Hence by the Calabi--Yau theorem for complex Monge--Amp\`{e}re equations there exists a smooth metric $\varphi_{\varepsilon}$ on $L$ such that 
\begin{equation}\label{e34}
i\Theta_{L,\varphi_{\varepsilon}}+\varepsilon\omega>0\quad\textrm{and}\quad(i\Theta_{L,\varphi_{\varepsilon}}+\varepsilon\omega)^{n}=C_{\varepsilon}\omega^{n}.
\end{equation}
Here $C_{\varepsilon}>0$ is a constant such that
\[
C_{\varepsilon}=\frac{\int_{X}(c_{1}(L)+\varepsilon\omega)^{n}}{\int_{X}\omega^{n}}\geqslant C\varepsilon^{n-\textrm{nd}(L)}
\]
for a universal constant $C$. Now let $\psi$ be a singular metric on $L$ with positive curvature current, which will be specified later. Consider the metric
\[
\phi_{\varepsilon}=\delta\varphi_{\varepsilon}+(1-\delta)\psi
\]
on $L$, where $\delta>0$ is a sufficiently small number which will be fixed later. Moreover, by construction,
\[
\begin{split}
i\Theta_{L,\phi_{\varepsilon}}+2\varepsilon\omega&=\delta(i\Theta_{L,\varphi_{\varepsilon}}+\varepsilon\omega)+(1-\delta)(i\Theta_{L,\psi}+\varepsilon\omega)+\varepsilon\omega \\ 
&\geqslant\delta(i\Theta_{L,\varphi_{\varepsilon}}+\varepsilon\omega)+\varepsilon\omega.
\end{split}
\]

Denote by $0<a_{1}\leqslant\cdots\leqslant a_{n}$ and $0<\hat{a}_{1}\leqslant\cdots\leqslant\hat{a}_{n}$, respectively, the eigenvalues of the curvature forms $i\Theta_{L,\varphi_{\varepsilon}}+\varepsilon\omega$ and $i\Theta_{L,\phi_{\varepsilon}}+2\varepsilon\omega$ at every point $x\in X$, with respect to the base K\"{a}hler metric $\omega(x)$. We find $\hat{a}_{j}\geqslant\delta a_{j}+\varepsilon$. On the other hand the Monge--Amp\`{e}re equation (\ref{e34}) tells us that
\begin{equation}\label{e35}
a_{1}\cdots a_{n}=C_{\varepsilon}\geqslant C\varepsilon^{n-\textrm{nd}(L)}
\end{equation}
everywhere on $X$.

Let us take a cohomology class $[\beta]\in H^{q}(X,K_{X}\otimes L\otimes\mathcal{E}(h))$. By using the de Rham--Weil isomorphism, we obtain a representative $\beta$ which is a smooth $L\otimes E$-valued $(n,q)$-form with bounded $L^{2}$-norm against $h$. Observe that
\[
\begin{split}
\|\beta\|^{2}_{\phi_{\varepsilon},h,\omega}&=\int_{X}|\beta|^{2}_{h,\omega}e^{-\delta\varphi_{\varepsilon}}e^{-(1-\delta)\psi}dV_{\omega} \\
&\leqslant(\int_{X}e^{-s\delta\varphi_{\varepsilon}}dV_{\omega})^{1/s}(\int_{X}|\beta|^{2t}_{h,\omega}dV_{\omega})^{1/t}(\int_{X}e^{-m(1-\delta)\psi}dV_{\omega})^{1/m}
\end{split}
\]
for any $s,t,m>1$ with $\frac{1}{s}+\frac{1}{t}+\frac{1}{m}=1$ by H\"{o}lder's inequality. We can take $t=1+\xi$ with $\xi$ small enough such that $\int_{X}|\beta|^{2t}_{h,\omega}dV_{\omega}<\infty$ by Corollary \ref{c21}. Then for an arbitrary corresponding pair $(s,m)$, due to the same reason as in (1), one can pick $\delta$ small enough to make $\int_{X}e^{-s\delta\varphi_{\varepsilon}}$ uniformly bounded. Whereas for any $x\in X$ there exists a suitable $\psi_{x}$ such that $e^{-m(1-\delta)\psi_{x}}$ is integrable around $x$ by the assumption that $L$ is analytic almost base point free. Then the compactness of $X$ implies a uniform $\psi$ such that $\int_{X}e^{-m(1-\delta)\psi}dV_{\omega}<\infty$. As a result, $\|\beta\|^{2}_{h,\phi_{\varepsilon},\omega}$ is uniformly bounded. Let $\Theta=i\Theta_{L,\phi_{\varepsilon}}+2\varepsilon\omega$, then
\[
\int_{X}\langle B^{-1}_{\Theta,\omega}\beta,\beta\rangle_{\phi_{\varepsilon},h,\omega}dV_{\omega}\leqslant\int_{X}\frac{1}{\hat{a}_{1}+\cdots+\hat{a}_{q}}|\beta|^{2}_{\phi_{\varepsilon},h,\omega}dV_{\omega}\leqslant\frac{1}{q\varepsilon}\|\beta\|^{2}_{h,\phi_{\varepsilon},\omega}<\infty.
\]
Apply the strong $L^{2}$-estimate condition, we obtain elements $u_{\varepsilon},v_{\varepsilon}$ and positive constant $C$ satisfying $\beta=\bar{\partial}u_{\varepsilon}+v_{\varepsilon}$ and
\[
\|u_{\varepsilon}\|^{2}_{\phi_{\varepsilon},h,\omega}+\frac{C}{\varepsilon}\|v_{\varepsilon}\|^{2}_{\phi_{\varepsilon},h,\omega}\leqslant\int_{X}\langle B^{-1}_{\Theta,\omega}\beta,\beta\rangle_{\phi_{\varepsilon},h,\omega}dV_{\omega}.
\]

Then a similar argument as in (1) implies that
\[
\begin{split}
\|v_{\varepsilon}\|^{2}_{\phi_{\varepsilon},h,\omega}&\leqslant\frac{\varepsilon}{C}\int_{X}\frac{1}{\hat{a}_{1}+\cdots+\hat{a}_{q}}|\beta|^{2}_{\phi_{\varepsilon},h,\omega}dV_{\omega} \\
&=\int_{X}\gamma_{\varepsilon}|\beta|^{2}_{h,\omega}e^{-(\delta\varphi_{\varepsilon}+(1-\delta)\psi)}dV_{\omega} \\
&\leqslant(\int_{X}e^{-s\delta\varphi_{\varepsilon}}dV_{\omega})^{1/s}(\int_{X}\gamma^{t}_{\varepsilon}|\beta|^{2t}_{h,\omega}dV_{\omega})^{1/t}(\int_{X}e^{-m(1-\delta)\psi}dV_{\omega})^{1/m}
\end{split}
\]
for any $s,t,m>1$ with $\frac{1}{s}+\frac{1}{t}+\frac{1}{m}=1$ by H\"{o}lder's inequality. Here 
\[
\gamma_{\varepsilon}=\frac{\varepsilon}{C(\hat{a}_{1}+\cdots+\hat{a}_{q})}\leqslant\min(1,C^{\prime}\varepsilon^{1-(n-\textrm{nd}(L))/q}(a_{q+1}\cdots a_{n})^{1/q}),
\]
hence converges almost everywhere to zero as $\varepsilon$ tends to zero whenever $q>n-\textrm{nd}(L)$. Remember that we have picked a collection of $(s,t,m,\delta;\psi)$ to make $\int_{X}e^{-s\delta\varphi_{\varepsilon}}dV_{\omega}$ uniformly bounded, and both of $\int_{X}|\beta|^{2t}_{h,\omega}dV_{\omega}$ and $\int_{X}e^{-m(1-\delta)\psi}dV_{\omega}$ bounded. As $\gamma_{\varepsilon}\leqslant1$, the Lebesgue dominated convergence theorem shows that
\[
\int_{X}\gamma^{t}_{\varepsilon}|\beta|^{2t}_{h,\omega}dV_{\omega}
\]
as well as $v_{\varepsilon}$ converges to zero as $\varepsilon$ tends to zero. 

In summary, $[\beta]=0$. The proof is complete.
\end{proof}
\end{theorem}

\section{Main theorem}
\label{sec:main} 

\subsection{Technical preparation}
Firstly we give a simple variant of \cite{Mat22}, Theorem 1.7 which is applied in the proof of our Theorem \ref{t12}. 

\begin{theorem}\label{t41}
Let $f:X\rightarrow Y$ be a K\"{a}hler fibration from a complex manifold $X$ to an analytic space $Y$. Let $L$ be a holomorphic line bundle on $X$ with $\kappa(L,f)\geqslant0$, and let $(E,h)$ be a smooth Hermitian vector bundle over $X$ that is Nakano semi-positive. For a (non-zero) section $s$ of some multiple $L^{m-1}$, the multiplication map induced by the tensor product with $s$
\begin{equation*}
\Phi:R^{q}f_{\ast}(K_{X}\otimes L\otimes E\otimes\mathscr{I}(f,\|L\|))\rightarrow R^{q}f_{\ast}(K_{X}\otimes L^{m}\otimes E\otimes\mathscr{I}(f,\|L^{m}\|))
\end{equation*}
is well-defined and injective for any $q\geqslant0$. In particular, $R^{q}f_{\ast}(K_{X}\otimes L\otimes E\otimes\mathscr{I}(f,\|L\|))$ is torsion-free for every $q$.
\begin{proof}
We only give a brief explanation that how to apply \cite{Mat22}. 

Firstly, since it only involves the local property, we can reduce to the setting that $f:X\rightarrow\Delta$ is a fibration from a K\"{a}hler manifold $X$ to a Stein subvariety $\Delta$ of $\mathbb{C}^{n}$. At this time, $\mathscr{I}(f,\|L\|)=\mathscr{I}(\|L\|)$. Let $\varphi$ be the associated metric, and it is enough to show the multiplication map induced by the tensor product with $s$
\[
H^{q}(X,K_{X}\otimes L\otimes E\otimes\mathscr{I}(\varphi))\rightarrow H^{q}(X,K_{X}\otimes L^{m}\otimes E\otimes\mathscr{I}(m\varphi))
\]
is injective. 

It furthermore involves a vector bundle $(E,h)$ here comparing with \cite{Mat22}. However it brings no essential trouble, since $(E,h)$ is supposed to be Nakano semi-positive. In particular, the key ingredients, such as Lemma 3.8, Propositions 3.15 and 3.17 in \cite{Mat22}, are still valid. The other things also extend without difficulty.
\end{proof}
\end{theorem}

Unfortunately, when $(E,h)$ is further a singular Nakano semi-positive vector bundle, the situation would be much more complicated. Hence we will apply the induction on dimension to show Theorem \ref{t14}. Whereas for Theorem \ref{t13}, we introduce the following lemma.

\begin{lemma}\label{l41}
Let $f:X\rightarrow Y$ be a fibration from a compact complex manifold $X$ to an irreducible projective variety $Y$, and let $A$ be an ample divisor on $Y$. Fix an integer $i$. Suppose that $\mathcal{F}$ is a coherent sheaf on $X$ with the property that
\[
H^{q}(X,\mathcal{F}\otimes\mathcal{O}_{X}(f^{\ast}mA))=0\quad\textrm{for all }q>i\textrm{ and  all }m\gg0. 
\]
Then $R^{q}f_{\ast}\mathcal{F}=0$ for every $q>i$.
\begin{proof}
Choose $m$ sufficiently large so that
\[
H^{p}(Y,R^{q}f_{\ast}\mathcal{F}\otimes\mathcal{O}_{Y}(mA))=0
\]
for all $p>0$ and $q\geqslant0$, and so that in addition the sheaves $R^{q}f_{\ast}\mathcal{F}\otimes\mathcal{O}_{Y}(mA)$, if non-zero, are globally generated. Then the Leray spectral sequence degenerates and shows that 
\[
H^{q}(X,\mathcal{F}\otimes\mathcal{O}_{Y}(f^{\ast}mA))=H^{0}(Y,R^{q}f_{\ast}\mathcal{F}\otimes\mathcal{O}_{Y}(mA)).
\]
We have arranged that the group on the right is non-zero if $R^{q}f_{\ast}\mathcal{F}\neq0$ and $m\gg0$ is sufficiently positive. But if $q>i$ this would violate the hypothesis, and the lemma follows.
\end{proof}
\end{lemma}

\subsection{Proof of the main theorem}
Now we should prove Theorems \ref{t12}, \ref{t13} and \ref{t14}.

\begin{theorem}[(=Theorem \ref{t12})]\label{t42}
Let $f:X\rightarrow Y$ be a K\"{a}hler fibration from a complex manifold $X$ to an analytic space $Y$. Suppose that the dimension of the general fiber is $l$. Let $(E,h)$ be a smooth Hermitian vector bundle that is Nakano semi-positive, and let $L$ be a holomorphic line bundle with $\kappa(L,f)\geqslant0$ on $X$. Then
\begin{equation*}
R^{q}f_{\ast}(K_{X}\otimes L\otimes E\otimes\mathscr{I}(f,\|L\|))=0
\end{equation*}
for $q>l-\kappa(L,f)$.   
\begin{proof}
Let $Z$ be the set of the critical value of $f$, which is a closed subvariety of $Y$. As a result, $X_{y}$ is a compact K\"{a}hler manifold when $y\in Y\setminus Z$. Then we apply Theorem \ref{t31} to obtain that 
\begin{equation*}
H^{q}(X_{y},K_{X_{y}}\otimes L|_{X_{y}}\otimes E|_{X_{y}}\otimes\mathscr{I}(\|L|_{X_{y}}\|))=0
\end{equation*}
for $q>l-\kappa(L|_{X_{y}})$.  

Now we claim that there exists a subset $V$ of $Y$ with zero Lebesgue measure, such that $\kappa(L|_{X_{y}})=\kappa(L,f)$ and
\begin{equation*}
\mathscr{I}(\|L|_{X_{y}}\|)=\mathscr{I}(f,\|L\|)|_{X_{y}}
\end{equation*}
when $y\in Y\setminus V$. In fact, the base-change theorem (c.f. \cite{Har77}, Theorem III.12.11) implies that for all $m$, there exists a closed subvariety $Z_{m}$ such that $f_{\ast}(L^{m})$ is locally free and
\[
f_{\ast}(L^{m})_{y}=H^{0}(X_{y},L^{m}|_{X_{y}})
\]
when $y\in Y\setminus Z_{m}$. Let $V=(\cup^{\infty}_{m=1}Z_{m})\cup Z$, which has zero Lebesgue measure. Then every section in $H^{0}(X_{y},L^{m}|_{X_{y}})$ extends locally and $\kappa(L|_{X_{y}})=\kappa(L,f)$ when $y\in Y\setminus V$. Hence 
\[
\mathscr{I}(\|L|_{X_{y}}\|)=\mathscr{I}(f,\|L\|)
\] 
by definition. 

Now on $Y\setminus V$ we obtain
\begin{equation*}
R^{q}f_{\ast}(K_{X}\otimes L\otimes E\otimes\mathscr{I}(f,\|L\|))=0
\end{equation*}
for $q>l-\kappa(L,f)$. We then conclude that this vanishing actually holds on the whole $Y$ since $R^{q}f_{\ast}(K_{X}\otimes L\otimes E\otimes\mathscr{I}(f,\|L\|))$ is torsion-free by Theorem \ref{t41}.
\end{proof}
\end{theorem}

\begin{theorem}[(=Theorem \ref{t13})]
Let $f:X\rightarrow Y$ be a fibration from an $n$-dimensional compact K\"{a}hler manifold $X$ to a projective manifold $Y$. Suppose that the dimension of the general fiber is $l$. Let $(E,h)$ be a singular Hermitian vector bundle that is strong Nakano semi-positive, and let $L$ be an analytic almost base point free line bundle on $X$. Then
\begin{equation*}
R^{q}f_{\ast}(K_{X}\otimes L\otimes\mathcal{E}(h))=0
\end{equation*}
for $q>l-\kappa(L,f)$.
\begin{proof}
For an arbitrary ample divisor $A$ on $Y$, $L\otimes\mathcal{O}_{X}(f^{\ast}A)$ is certainly analytic almost base point free, and we claim that for any $m\gg0$, 
\[
\textrm{nd}(L\otimes\mathcal{O}_{X}(f^{\ast}mA))\geqslant\kappa(L\otimes\mathcal{O}_{X}(f^{\ast}mA))\geqslant\kappa(L,f)+\dim Y.
\]
By Theorem \ref{t31},
\[
H^{q}(X,K_{X}\otimes L\otimes\mathcal{E}(h)\otimes\mathcal{O}_{X}(f^{\ast}mA))=0
\]
for $q>n-\dim Y-\kappa(L,f)=l-\kappa(L,f)$. Then we obtain the desired vanishing result due to Lemma \ref{l41}.

It remains to prove the claim. Indeed, for any $k$ such that $\textrm{rank}\,f_{\ast}(L^{k})\sim k^{\kappa(L,f)}$, by Serre's vanishing theorem we have
\[
h^{0}(X,L^{k}\otimes\mathcal{O}_{X}(f^{\ast}kmA))=h^{0}(Y,f_{\ast}(L^{k})\otimes\mathcal{O}_{Y}(kmA))=\chi(Y,f_{\ast}(L^{k})\otimes\mathcal{O}_{Y}(kmA))
\]
provided $m\gg0$. Then due to the Hirzebruch--Riemann--Roch formula we obtain that 
\[
h^{0}(X,L^{k}\otimes\mathcal{O}_{X}(f^{\ast}kmA))\sim\textrm{rank}\,f_{\ast}(L^{k})\cdot k^{\dim Y}\sim k^{\kappa(L,f)+\dim Y},
\]
which implies the desired result.
\end{proof}
\end{theorem}

\begin{theorem}[(=Theorem \ref{t14})]
Let $f:X\rightarrow Y$ be a fibration from an $n$-dimensional projective manifold $X$ to an analytic space $Y$. Suppose that the dimension of the general fiber is $l$. Let $(E,h)$ be a singular Hermitian vector bundle that is strong Nakano semi-positive, and let $L$ be an $f$-nef line bundle on $X$. Then
\begin{equation*}
R^{q}f_{\ast}(K_{X}\otimes L\otimes\mathcal{E}(h))=0
\end{equation*}
for $q>l-\kappa(L,f)$.
\begin{proof}
Let $A$ be a very ample line bundle on $X$. Let $H_{1}, \ldots, H_{n}$ be generic hypersurfaces in $|A|$, and let
\[
S_{r}:=H_{1} \cap \cdots \cap H_{r}
\]
with the convention that \(S_{0}=X\). Then there exists a natural surjective morphism \(f_{r}: S_{r} \to Y_{r}:=f(S_{r})\).

By the standard exact sequence
\[
0\to K_{X}\otimes L\otimes\mathcal{E}(h)\to K_{X}\otimes L\otimes A\otimes\mathcal{E}(h)\to K_{S_{1}}\otimes L|_{S_{1}}\otimes \mathcal{E}(h)|_{S_{1}} \to 0,\]
we obtain the long exact sequence
\[
\begin{aligned}
\cdots \to & R^{q-1}\left(f_{1}\right)_{*}\left(\left.K_{S_{1}} \otimes L\right|_{S_{1}} \otimes \mathcal{E}(h)|_{S_{1}}\right) \to R^{q} f_{*}\left(K_{X} \otimes L \otimes \mathcal{E}(h)\right) \to \\
& R^{q} f_{*}\left(K_{X} \otimes L \otimes A \otimes \mathcal{E}(h)\right) \to R^{q}\left(f_{1}\right)_{*}\left(\left.K_{S_{1}} \otimes L\right|_{S_{1}} \otimes \mathcal{E}(h)|_{S_{1}}\right) \to \cdots.
\end{aligned}
\]

Observe that we can pick \(A\), so that \(R^{q} f_{*}(K_{X} \otimes L \otimes A \otimes \mathcal{E}(h)) = 0\) for \(q>0\) by Serre's vanishing theorem. Hence
\[
R^{q} f_{*}\left(K_{X} \otimes L \otimes \mathcal{E}(h)\right)=R^{q-1}\left(f_{1}\right)_{*}\left(K_{S_{1}} \otimes L|_{S_{1}} \otimes \mathcal{E}(h)|_{S_{1}}\right) .
\]
Repeating this procedure, we will eventually obtain that
\[
R^{q} f_{*}\left(K_{X} \otimes L \otimes \mathcal{E}(h)\right)=R^{1}\left(f_{q-1}\right)_{*}\left(\left.K_{S_{q-1}} \otimes L\right|_{S_{q-1}} \otimes \mathcal{E}(h)|_{S_{q-1}}\right) .
\]

Since the relative Iitaka dimension is non-decreasing in restriction, \(L|_{S_{q-1}}\) is actually \((f_{q-1})\)-big when \(q>l-\kappa(L, f)\). In this situation we will prove that
\[
\mathcal{F}:=R^{1}\left(f_{q-1}\right)_{*}\left(\left.K_{S_{q-1}} \otimes L\right|_{S_{q-1}} \otimes \mathcal{E}(h)|_{S_{q-1}}\right)=0,
\]
which leads to the desired vanishing result. By definition, $\mathcal{F}$ is associated with the presheaf
\[
U\mapsto H^{1}(f^{-1}_{q-1}(U),K_{S_{q-1}}\otimes L|_{S_{q-1}}\otimes\mathcal{E}(h)|_{S_{q-1}})
\]
for any Stein open set $U$. Fix a K\"{a}hler metric $\omega$ on $f^{-1}_{q-1}(U)$, and take an arbitrary 
\[
[\beta]\in H^{1}(f^{-1}_{q-1}(U),K_{S_{q-1}}\otimes L|_{S_{q-1}}\otimes\mathcal{E}(h)|_{S_{q-1}}),
\]
which is represented by an $L\otimes E$-valued $\bar{\partial}$-closed $(n-q+1,1)$-form $\beta$ satisfying
\[
\int_{f^{-1}_{q-1}(U)}|\beta|^{2}_{h,\omega}dV_{\omega}<\infty.
\]
Then for any singular metric $\varphi$ on $L$,
\[
\int_{f^{-1}_{q-1}(U)}|\beta|^{2}_{\varphi,h,\omega}dV_{\omega}\leqslant(\int_{f^{-1}_{q-1}(U)}e^{-s\varphi}dV_{\omega})^{1/s}(\int_{f^{-1}_{q-1}(U)}|\beta|^{2t}_{h,\omega}dV_{\omega})^{1/t}
\]
for any $s,t>1$ with $\frac{1}{s}+\frac{1}{t}=1$ by H\"{o}lder's inequality. Pick $t=1+\xi$ with $\xi$ small enough such that $\int_{f^{-1}_{q-1}(U)}|\beta|^{2t}_{h,\omega}dV_{\omega}<\infty$. On the other hand, since $L$ is nef and big on $f^{-1}_{q-1}(U)$, by \cite{Dem12} there exists a singular metric $\varphi_{0}$ such that $i\Theta_{L,\varphi_{0}}>\delta\omega$ for some positive $\delta$ and
\[
\int_{f^{-1}_{q-1}(U)}e^{-s\varphi_{0}}dV_{\omega}<\infty
\]
for the corresponding $s$. In summary,
\[
\int_{f^{-1}_{q-1}(U)}|\beta|^{2}_{\varphi_{0},h,\omega}dV_{\omega}<\infty.
\] 

Take $\Theta:=i\Theta_{L,\varphi_{0}}$, and let $\delta<a_{1}\leqslant\cdots\leqslant a_{n}\leqslant+\infty$ be the eigenvalues of $\Theta$ at every point $x\in f^{-1}_{q-1}(U)$ with respect to $\omega(x)$. Then
\[
\int_{f^{-1}_{q-1}(U)}\langle B^{-1}_{\Theta,\omega}\beta,\beta\rangle_{\varphi_{0},h,\omega}dV_{\omega}\leqslant\frac{1}{\delta}\int_{f^{-1}_{q-1}(U)}|\beta|^{2}_{\varphi_{0},h,\omega}dV_{\omega}<\infty.
\]
By assumption $(E,h)$ satisfies the strong $L^{2}$-estimate condition, so there exist elements $u_{\varepsilon},v_{\varepsilon}$ and positive constant $C$ with $\beta=\bar{\partial}u_{\varepsilon}+v_{\varepsilon}$ and
\[
\|u_{\varepsilon}\|^{2}_{\varphi_{0},h,\omega}+\frac{C}{\varepsilon}\|v_{\varepsilon}\|^{2}_{\varphi_{0},h,\omega}\leqslant\int_{f^{-1}_{q-1}(U)}\langle B^{-1}_{\Theta,\omega}\beta,\beta\rangle_{\varphi_{0},h,\omega}dV_{\omega}.
\]
Note for every $\varepsilon$ we take $\varphi_{0}$ as the corresponding metric here. Now 
\[
\|v_{\varepsilon}\|^{2}_{\varphi_{0},h,\omega}\leqslant\frac{\varepsilon}{C}\int_{f^{-1}_{q-1}(U)}\langle B^{-1}_{\Theta,\omega}\beta,\beta\rangle_{\varphi_{0},h,\omega}dV_{\omega},
\]
which converges to zero as $\varepsilon$ tends to zero. In conclusion, $[\beta]=0$. The proof is complete.
\end{proof}
\end{theorem}

\end{document}